\documentclass{article}

\usepackage{amsmath}
\usepackage{amsfonts}

\usepackage{listings}
\usepackage[numbered]{mcode} 
  
   \usepackage{longtable}
   \usepackage{hyperref}
  \usepackage{pgfplots}
  
  \usepackage{multirow}

  \DeclareMathOperator*{\argmin}{\arg \min}%
  \DeclareMathOperator*{\argmax}{\arg \max}%
  \newcommand*\cpp{C\kern-0.2ex\raisebox{0.4ex}{\scalebox{0.8}{+\kern-0.4ex+}}}
  
\title{A \textbf{\underline{Flex}}ible Primal-Dual Tool\textbf{\underline{Box}} \\\large Technical Report}
\author{Hendrik Dirks\footnotemark[1]}

\begin{document}
	
\maketitle

\footnotetext[1]{Institute for Computational and Applied Mathematics and Cells in Motion Cluster of Excellence, University of M{\"u}nster, Orl\'{e}ans-Ring 10, 48149 M{\"u}nster, Germany, Email: \href{mailto:hendrik.dirks@wwu.de}{hendrik.dirks@wwu.de}}

\begin{abstract}
	\textbf{FlexBox} is a flexible MATLAB toolbox for finite dimensional convex variational problems in image processing and beyond. Such problems often consist of non-differentiable parts and involve linear operators. The toolbox uses a primal-dual scheme to avoid (computationally) inefficient operator inversion and to get reliable error estimates. From the user-side, \textbf{FlexBox} expects the primal formulation of the problem, automatically decouples operators and dualizes the problem. For large-scale problems, \textbf{FlexBox} also comes with a \cpp-module, which can be used stand-alone or together with MATLAB via MEX-interfaces. Besides various pre-implemented data-fidelities and regularization-terms, \textbf{FlexBox} is able to handle arbitrary operators while being easily extendable, due to its object-oriented design.\\
	The toolbox is available at \href{http://www.flexbox.im}{http://www.flexbox.im}
\end{abstract}

\section{Introduction}
Many variational problems in image processing can be written in the form
\begin{align}
	\argmin_{x} G(x) + F(Ax),
	\label{generalPrimalFormulation}
\end{align}
where $A$ denotes some linear operator (see \ref{appendix}) and both $G$ and $F$ are proper, convex and lower-semicontinuous functions. Problem \eqref{generalPrimalFormulation} refers to the so-called \textit{primal} formulation of the minimization problem and $x$ is known as the primal variable. It can be shown (see \cite{rockafellar2015convex}) that minimizing \eqref{generalPrimalFormulation} is equivalent to solving the \textit{primal-dual} or \textit{saddle-point} formulation
\begin{align}
	\argmin_{x}\argmax_{y} G(x) + \langle y,Ax\rangle - F^*(y).
\label{generalPrimalDualFormulation}
\end{align}
Here, $F^*$ refers to the convex conjugate (see \ref{appendix}) of $F$ and $y$ is denoted as the dual variable. It the recent years, algorithms like ADMM \cite{parikh2014proximal,wahlberg2012admm} or primal-dual \cite{pock2009algorithm,chambolle2011first,zhang2011unified,esser2010general} for efficiently solving these saddle-point problems have become very popular. \textbf{FlexBox} makes use of the latter, which can be sketched up as follows:\\
For $\tau,\sigma>0$ and a pair $(\hat{x}^0,y^0)\in\mathcal{X}\times\mathcal{Y}$ we iteratively solve:
\begin{align}
	y^{k+1} &= prox_{\sigma F^*}(y^k+\sigma A\hat{x}^k)\\
	x^{k+1} &=  prox_{\tau G}(x^k-\tau A^T y^{k+1})\\
	\hat{x}^{k+1} & = 2 x^{k+1} -x^k
	\label{pdAlgorithm}
\end{align}
Here, $prox_{\tau G}$ (resp. $\sigma F^*$) denotes the \textit{proximal} or \textit{resolvent} operator
\begin{align*}
	prox_{\tau G}(y) = \left( I+\tau \partial G\right)^{-1}(y) := \argmin_v \left\{ \frac{\left\| v-y \right\|_2^2}{2}+\tau G(v) \right\} ,
\end{align*}
which can be interpreted as a compromise between minimizing $G$ and being close to the input argument $y$. The efficiency of primal-dual algorithms relies on the fact that the \textit{prox-problems} are computationally efficient to solve.
\subsection{Contribution}
Since primal-dual algorithms have been extensively applied to all classes of convex optimization problems, we found that people are spending a lot of effort on calculating convex conjugates or solutions for prox-problems repeatingly for similar problems. Let us consider, for example, the isotropic total variation $\|\nabla u\|_{1,2}$, where the convex conjugate is an indicator function of the L$^\infty$-ball and the solution of the prox-problem is a point-wise projection onto L$^2$-balls. These results hold not only for $A=\nabla$, but for arbitrary operators. \textbf{FlexBox} makes use of this generalization and simply works on the level of terms in the primal problem. After adding a certain term, \textbf{FlexBox} automatically decouples operators, creates dual variables and calculates step-sizes ($\tau, \sigma$). \textbf{FlexBox} already contains a variety of data-fidelity (e.g. L$^1$, L$^2$, Kullback-Leibler) and regularization terms (e.g. L$^2$, TV, Laplace, curl), but is also compatible with user-defined operators. Moreover, the class-based structure allows easy extension and creation of custom terms. A full list of available terms can be found in Table \ref{tab:listPrimalDualNotation}.\\
Core components of \textbf{FlexBox} are written in MATLAB, but there exists an optional \cpp-module to improve compatibility and runtime. This module can be compiled and will afterwards be used automatically via a MEX-interface. The \cpp-module can also be used without MATLAB (but e.g. with OpenCV \cite{bradski2000opencv,bradski2008learning}), but does not have the full functionality and whole variety of terms included.

\section{Architecture and Features}
\paragraph{Basic Idea:}
To work with \textbf{FlexBox}, the user has to derive the primal formulation of the variational problem first. Afterwards, each primal variable including its dimensions is added to \textbf{FlexBox}. Then, the primal problem is put into the toolbox term by term, using the implemented functional terms from Table \ref{tab:listPrimalDualNotation}. Once the computation is finished, the result stored in the primal variables can be requested and used.

\paragraph{Design:} 
\textbf{FlexBox} is designed as one core class, which holds a list of functional terms. For each term the toolbox decides whether it has to be dualized and creates necessary dual variables. The main object holds the data of primal and dual variables $x_i$ and $y_i$ and all parameters. Terms always correspond to at least one primal variable, whereas dual terms also correspond to at least one dual variable and contain the involved operators.\\ 
In the primal-dual algorithm \eqref{pdAlgorithm}, applications of the operator can be decoupled defining 
$$\tilde{y} := y^k+\sigma A\hat{x}^k, \text{ and } \tilde{x} := x^k-\tau A^T y^{k+1}.$$
Since dual terms bind the operators, calculations of $\tilde{y}$ and $\tilde{x}$ are done inside the corresponding dual terms, accessing the variables held by the core class. Finally, primal and dual terms specify prox-methods to solve the arising prox-problems while again accessing the variables held by the \textbf{FlexBox} core class.
\paragraph{Parameters:}
To ensure convergence, the parameters $\tau$ and $\sigma$ have to fulfill $\tau\sigma\|A\| < 1$. A static choice of these parameters might lead to slow convergence speed, because it is dominated by the smallest possible value along all operators. A fully automatic strategy for finding \textit{optimal} parameters can be extrapolated from \cite{pock2011diagonal}. Summing up the absolute values of row elements (for $\sigma_i$) resp. column elements (for $\tau_i$) leads to custom parameters for each term in the functional. This strategy is inherited by \textbf{FlexBox}. 
\paragraph{Stopping Criterion:}
As a powerful stopping criterion \textbf{FlexBox} uses the primal-dual residual, proposed by Goldstein, Esser and Baraniuk \cite{goldstein2013adaptive}, which can be calculated after the (k+1)-th iteration as
\begin{align*}
	p^k := \left|\frac{x^k-x^{k+1}}{\tau} - A^T(y^k-y^{k+1}) \right|,\quad
	d^k := \left|\frac{y^k-y^{k+1}}{\sigma} - A(x^k-x^{k+1}) \right|,
\end{align*}
with $\left|\cdot\right|$ being the sum of absolute values. We denote $p^k$ as the primal and $d^k$ as the dual residual. The total residual is then given by the sum of primal- and dual residual, which is afterwards scaled with the size of the problem and number of variables. Since evaluating the residual is computationally expensive, it is regularly computed after a fixed number of iterations (default 100).\\
Besides evaluating the primal-dual residual, \textbf{FlexBox} automatically stops after a static number of iterations (default 10000) and can be continued afterwards.
\section{Examples}

\subsection{Rudin-Osher-Fatemi}
The Rudin-Osher-Fatemi model \cite{rudin1992nonlinear} has very popular applications in image denoising. The primal formulation reads
\begin{align}
	\argmin_u \frac{1}{2}\|u-f\|_2^2 + \alpha \|\nabla u\|_{1,2},
\end{align}
where the first part fits the unknown $u$ to the given input image $f$ and the second part refers to the isotropic total variation, which penalizes the total sum of jumps in the solution. Minimizing this problem with \textbf{FlexBox} can be done with the following lines of code:
\vspace{-4mm}
\begin{lstlisting} 
%Begin: Code example
main = flexBox;

numberU = main.addPrimalVar(size(f));

main.addTerm(L2dataTerm(1,imageNoisy),numberU);
main.addTerm(L1gradientIso(0.08,size(f)),numberU);

main.runAlgorithm;

result = main.getPrimal(numberU);
%End: Code example
\end{lstlisting}
To keep this example short, we omitted parts of the code where the image $f$ is read. Let us begin in line $2$, which initializes a \textbf{FlexBox} object and saves it to the variable $main$. Line 4 then adds the primal objective variable $u$, which has the same size as the input image $f$. The toolbox returns the internal number of this primal variable, which is saved in $numberU$. \\
In line 6 and 7, the L$^2$-data-term with weight $1$ (the weight is divided by 2 internally for L$^2$ terms) and corresponding image $f$ is added. Moreover, the isotropic TV-term with weight $0.08$ is pushed into the framework. Note that the function $addTerm$ always requires a functional part and the internal number of the corresponding primal variable(s). \\
The function call in line $9$ finally starts the calculation. Once this is finished we transfer the solution into the variable $result$ in line $11$. 

\subsection{Optical Flow}
Estimating the motion between two consecutive images $f_1$ and $f_2$ based on the displacement of intensities in both images is called \textit{optical-flow} estimation. The unknown velocity field $\boldsymbol{v}$ is usually connected to the image by the brightness-constancy-assumption $f_2(x+\boldsymbol{v})-f_1(x)=0$. This formulation is non-linear in terms of $\boldsymbol{v}$ and therefore linearized (see e.g. \cite{zach2007duality,dirks}). A corresponding variational problem incorporating total variation regularization can be written as
\begin{align}
	\argmin_{\boldsymbol{v}=(v_1,v_2)} \frac{1}{2}\|f_2-f_1 + \nabla f_2\cdot\boldsymbol{v}\|_1 + \alpha_1 \|\nabla v_1\|_{1,2} + \alpha_2\|\nabla v_2\|_{1,2}.
\end{align}
Solving this problem with \textbf{FlexBox} can be done in a similar manner as for the ROF example:
\vspace{-4mm}
\begin{lstlisting} 
%Begin: Code example
main = flexBox;

numberV1 = main.addPrimalVar(size(f1));
numberV2 = main.addPrimalVar(size(f2));

%add optical flow data term
main.addTerm(L1opticalFlowTerm(1,f1,f2),[numberV1,numberV2]);

%add regularizers - one for each component
main.addTerm(L1gradientIso(0.05,size(f1)),numberV1);
main.addTerm(L1gradientIso(0.05,size(f1)),numberV2);

main.runAlgorithm;

resultV1 = main.getPrimal(numberV1);
resultV2 = main.getPrimal(numberV2);
%End: Code example
\end{lstlisting}
In lines 4 and 5 primal variables for both components of the velocity fields are added. Afterwards, in lines 8, 11 and 12 the data term and regularizers for both components are inserted. Please note that the optical flow term now refers to two primal variables written as the vector $[numberV1,numberV2]$. Afterwards, the algorithm is started and both results are retrieved.

\subsection{Segmentation}
Dividing an image into different regions is called \textit{segmentation}. Assuming that the image $f$ consists of $k$ different regions, each of them having a mean intensity $c_i$ for $i=1,\ldots,k$, the segmentation problem including total variation regularization can be written as
 \begin{align}
 \argmin_{u=u_1,\ldots,u_k}& \sum_{i=1}^{k}u_i\frac{1}{2}\|f-c_i\|_2^2 + \alpha \|\nabla u_i\|_{1,2},\\
 &s.t.\quad u_i\geq 0,\quad \sum_{i=1}^{k} u_i=1,
 \end{align}
where $u$ is a labeling vector (see \cite{mumford1989optimal,zach2008fast}). This labeling formulation is a convex relaxation of the integer assignment $u_i\in\{0,1\}$. The MATLAB implementation of this problem is again rather short:
\vspace{-4mm}
\begin{lstlisting} 
%Begin: Code example
numberOfLabels = 3;
dims = size(image);
labels = rand(numberOfLabels,1);

main = flexBox;

for i=1:numberOfLabels
	main.addPrimalVar(size(image));
end

%init data term
main.addTerm(labelingTerm(1,image,labels),1:numberOfLabels);

for i=1:numberOfLabels
	main.addTerm(L1gradientIso(0.5,size(image)),i);
end

main.runAlgorithm;

for i=1:numberOfLabels
	labelMatrix(:,:,i) = main.getPrimal(i);
end
%End: Code example
\end{lstlisting}
We begin by choosing a fixed number of regions and choose the \textit{mean} intensity in each region as random values in line $4$. Afterwards, the main object is initialized and a primal variable for each $u_i$ is created in lines 8-10. In line 13, we create the labeling term using the previously defined labels and draw a connection for primal variables \textit{1:numberOfLabels}. The loop in line 13-15 generates a total variation regularizer for each of the primal variables. The problem is solved in line 19 and we save each labeling function as a layer in a 3d matrix.

\section{Features}

\subsection{General}
\begin{itemize}
	\item \textbf{FlexBox} can be stopped at any time and afterwards continued from the current state. Due to this feature, we are able to change parameters after the problem has converged and use the current state as an initial guess for the next run.
	\item The toolbox supports arbitrary user-defined operators consisting of blocks $A_i$ (see \textit{General operator regularization} in Table \ref{tab:listPrimalDualNotation} ) that can be submitted by defining $A$ as a cell-array of blocks, where each \textit{numPrimals} element correspond to one row in the overall operator. An example for a user-defined operator is provided with the software package.
\end{itemize}
	\newpage
\subsection{Available Terms}
 \vspace*{-5mm}
\begingroup
\tiny
\setlength{\LTleft}{-15cm plus -1fill}
\setlength{\LTright}{\LTleft}
\begin{longtable}{|c|c|c|ll|}
	\hline
	\multicolumn{5}{|c|}{\normalsize \textbf{List of implemented terms}} \\
	\hline
	\textbf{Term} & \textbf{\cpp} & \textbf{Classname} & \multicolumn{2}{|c|}{\textbf{Parameters}}\\
	\hline
	\multicolumn{5}{|c|}{\textbf{Data-fidelity}} \\
	\hline
	
	$\alpha\|u-f\|_1$ & \checkmark & L1dataTerm(alpha,f) & f:& input data\\
	&&&$\alpha$:& weight\\
	\hline
	$\alpha\|Au-f\|_1$ & \checkmark & L1dataTermOperator(alpha,A,f) & f:& input data\\
	&&&$\alpha$:& weight\\
	&&&A:& operator\\
	\hline
	$\frac{\alpha}{2}\|u-f\|_2^2$ & \checkmark & L2dataTerm(alpha,f) & f:& input data\\
	&&&$\alpha$:& weight\\
	\hline
	$\frac{\alpha}{2}\|Au-f\|_2^2$ & \checkmark & L2dataTermOperator(alpha,A,f) & f:& input data\\
	&&&$\alpha$:& weight\\
	&&&A:& operator\\
	\hline
	$\int Au-f+f\log\frac{f}{Au}$ & \checkmark & KLdataTermOperator(alpha,A,f) & f:& input data\\
	s.t. $u\geq0$&&&$\alpha$:& weight\\
	&&&A:& operator\\
	\hline
	$\alpha\|\nabla f_2\cdot \boldsymbol{v}+f_2-f_1\|_1$  & \checkmark & L1opticalFlowTerm(alpha,f$_1$,f$_2$) & f$_1$,f$_2$:& images\\
	&&&$\alpha$:& weight\\
	\hline
	$\frac{\alpha}{2}\|\nabla f_2\cdot \boldsymbol{v}+f_2-f_1\|_2^2$ & \checkmark & L2opticalFlowTerm(alpha,f$_1$,f$_2$) & f$_1$,f$_2$:& images\\
	&&&$\alpha$:& weight\\
	\hline
	$\sum_{i=1}^n \langle u_i,f_i \rangle$ & \texttimes & labelingTerm(alpha,f,l) & f:& image\\
	s.t. $u_i\geq 0,\sum u_i=1$&&&$\alpha$:& weight\\
	s.t. $f_i = (f-l_i)^2$&&&l:& label vector\\
	&&&dims:& dimensions of u\\
	\hline
	\multicolumn{5}{|c|}{\textbf{Gradient regularization}} \\
	\hline
	$\alpha\|\nabla u\|_{1,1}$ & \checkmark & L1gradientAniso(alpha,dims)&$\alpha$:& weight\\
	&&&dims:& dimensions of u\\
	\hline
	$\alpha\|\nabla u\|_{1,2}$ & \checkmark & L1gradientIso(alpha,dims)&$\alpha$:& weight\\
	&&&dims:& dimensions of u\\
	\hline
	$\frac{\alpha}{2}\|\nabla u\|_2^2$ & \checkmark & L2gradient(alpha,dims)&$\alpha$:& weight\\
	&&&dims:& dimensions of u\\
	\hline
	$\alpha\|\nabla u\|_{H_\epsilon}$ & \texttimes & huberGradient(alpha,dims,epsi)&$\alpha$:& weight\\
	&&&dims:& dimensions of u\\
	&&&epsi:& $\epsilon$ for Huber-norm\\
	\hline
	$\alpha\|\nabla u\|_{F}$ & \checkmark & frobeniusGradient(alpha,dims)&$\alpha$:& weight\\
	&&&dims:& dimensions of u\\
	\hline
	$\alpha\|\nabla (u-w)\|_{1,1}$ & \checkmark & L1gradientDiffAniso(alpha,dims)&$\alpha$:& weight\\
	&&&dims:& dimensions of u\\
	\hline
	$\alpha\|\nabla (u-w)\|_{1,2}$ & \checkmark & L1gradientDiffIso(alpha,dims)&$\alpha$:& weight\\
	&&&dims:& dimensions of u\\
	\hline
	$\alpha\|\nabla u-w\|_{1,1}$ & \checkmark & L1secondOrderGradientAniso(alpha,dims)&$\alpha$:& weight\\
	&&&dims:& dimensions of u\\
	\hline
	$\alpha\|\nabla u-w\|_{1,2}$ & \checkmark & L1secondOrderGradientIso(alpha,dims)&$\alpha$:& weight\\
	&&&dims:& dimensions of u\\
	\hline
	\multicolumn{5}{|c|}{\textbf{General operator regularization}} \\
	\hline
	$\alpha\|A u\|_{1,1}$ & \checkmark & L1operatorAniso(alpha,numPrimals,A)&$\alpha$:& weight\\
	&&& numPrimals:& \# corresp. primals\\
	&&& A:& operator(s)\\
	\hline
	$\alpha\|A u\|_{1,2}$ & \checkmark & L1operatorIso(alpha,numPrimals,A)&$\alpha$:& weight\\
	&&& numPrimals:& \# corresp. primals\\
	&&& A:& operator(s)\\
	\hline
	$\frac{\alpha}{2}\|A u\|_2^2$ & \checkmark & L2operator(alpha,numPrimals,A)&$\alpha$:& weight\\
	&&& numPrimals:& \# corresp. primals\\
	&&& A:& operator(s)\\
	\hline
	$\alpha\|Au\|_{F}$ & \checkmark & frobeniusOperator(alpha,numPrimals,A)&$\alpha$:& weight\\
	&&& numPrimals:& \# corresp. primals\\
	&&& A:& operator(s)\\
	\hline
	$\alpha \langle b,Au\rangle$ & \texttimes & innerProductOperator(alpha,A,b)&$\alpha$:& weight\\
	&&&A:& operator\\
	&&&b:& vector b\\
	\hline
	\multicolumn{5}{|c|}{\textbf{Vector-field regularization}} \\
	\hline
	$\alpha\|curl(\boldsymbol{v})\|_{1}$ & \checkmark & L1curl(alpha,dims)&$\alpha$:& weight\\
	&&&dims:& dimensions of e.g. $v_1$\\
	\hline
	$\frac{\alpha}{2}\|curl(\boldsymbol{v})\|_{2}^2$ & \checkmark & L2curl(alpha,dims)&$\alpha$:& weight\\
	&&&dims:& dimensions of e.g. $v_1$\\
	\hline
	$\alpha\|\nabla\cdot\boldsymbol{v}\|_{1}$ & \checkmark & L1divergence(alpha,dims)&$\alpha$:& weight\\
	&&&dims:& dimensions of e.g. $v_1$\\
	\hline
	$\frac{\alpha}{2}\|\nabla\cdot\boldsymbol{v}\|_{2}^2$ & \checkmark & L2divergence(alpha,dims)&$\alpha$:& weight\\
	&&&dims:& dimensions of e.g. $v_1$\\
	\hline
	\multicolumn{5}{|c|}{\textbf{Other regularization}} \\
	\hline
	$\alpha\|u\|_{1}$ & \checkmark & L1identity(alpha,dims)&$\alpha$:& weight\\
	&&&dims:& dimensions of u\\
	\hline
	$\frac{\alpha}{2}\|u\|_2^2$ & \checkmark & L2identity(alpha,dims)&$\alpha$:& weight\\
	&&&dims:& dimensions of u\\
	\hline
	$\alpha \langle b,\nabla u\rangle$ & \texttimes & innerProductGradient(alpha,dims,b)&$\alpha$:& weight\\
	&&&dims:& dimensions of u\\
	&&&b:& vector b\\
	\hline
	$\alpha \langle b,\nabla (u-w)\rangle$ & \texttimes & innerProductGradientDiff(alpha,dims,b)&$\alpha$:& weight\\
	&&&dims:& dimensions of u\\
	&&&b:& vector b\\
	\hline
	
	\caption{List of terms currently available in \textbf{FlexBox} (incomplete).}	
	\end{longtable}
\label{tab:listPrimalDualNotation}
	\endgroup

\begin{appendix}
	\section{General Considerations}
	\label{appendix}
	
	\begin{itemize}
		\item Throughout this report, we consider finite dimensional optimization problems with primal variables $x\in\mathcal{X}$ and $y\in\mathcal{Y}$.
		\item The sets $\mathcal{X}\subset\mathbb{R}^{N}$ and $\mathcal{Y}\subset\mathbb{R}^{M}$ are assumed to be convex.
		\item We do not explicitly distinguish between a variable $x_c$ on a regular cartesian grid $N_1\times\ldots\times N_d$ (s.t. $N_1\cdot\ldots\cdot N_d=N$) and its vectorized equivalent $x\in\mathbb{R}^{N} $, which is gained by concatenating $x$ along the first dimension.
		\item The scalar product of two vectors on $\mathcal{X}$ is defined as 
		\begin{align*}
			\langle u,v \rangle = \sum_{i} u_{i} v_i,\quad u,v\in \mathcal{X}.
		\end{align*}
		\item \textbf{FlexBox} expects the matrix representation of the linear operator $A$, hence $A\in\mathbb{R}^{M\times N}$ and $A^*=A^T$. Note that the evaluation of a linear operator can always be written as a matrix-vector multiplication by applying the Kronecker product (see e.g. \cite{steeb1991kronecker}).
		\item The function $F^*$ refers to the convex conjugate or Legendre-Fenchel transformation of $F$ and is defined as
		\begin{align*}
			F^*(y^*) := \sup_{y\in\mathcal{Y}} \langle y^*,y\rangle - F(y).
		\end{align*}
	\end{itemize}
	
\end{appendix}

\section*{Acknowledgements}
\begin{itemize}
	\item The work has been supported by ERC via Grant EU FP 7 - ERC Consolidator Grant 615216 LifeInverse.
	\item Thanks to Eva-Maria Brinkmann, Janic F{\"o}cke, Lena Frerking, Julian Rasch and Carolin Ro{\ss}manith for testing \textbf{FlexBox} and discussing useful features.
	\item \textbf{FlexBox} makes use of the primal-dual algorithm working on concrete functional terms. A MATLAB toolbox consisting of a more general class of solvers for optimization problems is OOMFIP, developed by \href{http://www.damtp.cam.ac.uk/people/mb941/}{Martin Benning}. Thank you for your advices!
\end{itemize}	

\newpage
\bibliographystyle{plain}
\bibliography{citations}

	
\end{document}